%% file: paper.tex
\documentclass[twocolumn,10pt]{article}

\usepackage{amsmath,algpseudocode,dpr_fec,hyperref}
\usepackage{fec_tikz}
\hypersetup{colorlinks,linkcolor=blue,urlcolor=blue,citecolor=blue}
\usepackage{amsfonts}
\usepackage{dsfont}
\usepackage{geometry}
\usepackage{algorithm}
\newtheorem{thm}{Theorem}[section]

\newtheorem{condition}[thm]{Condition}

\usepackage{titlesec}
\titleformat{\section}{\large\bfseries}{\thesection}{1em}{}
\titleformat{\subsection}{\normalsize\bfseries}{\thesubsection}{1em}{}

\newif\ifappendix
\appendixtrue

\begin{document}

\bcenter
  \textbf{\large Adaptive Stochastic Optimization}
  \smallskip
  
  Frank E.~Curtis $\circ$ Katya Scheinberg
  \smallskip
  
  \today
\ecenter

\paragraph{Abstract.}  Optimization lies at the heart of machine learning and signal processing.  Contemporary approaches based on the stochastic gradient method are non-adaptive in the sense that their implementation employs prescribed parameter values that need to be tuned for each application.  This article summarizes recent research and motivates future work on adaptive stochastic optimization methods, which have the potential to offer significant computational savings when training large-scale systems.

\section{Introduction}

The successes of stochastic optimization algorithms for solving problems arising machine learning (ML) and signal processing (SP) are now widely recognized.  Scores of articles have appeared in recent years as researchers aim to build upon fundamental methodologies such as the stochastic gradient method (SG) \cite{RobbMonr51}.  The motivation and scope of much of these efforts have been captured in various books and review articles; see, e.g.,  \cite{BottCurtNoce18,CurtSche17,GoodBengCour16,MohrRostTalw18,SraNowoWrig11}.

Despite these advances and accumulation of knowledge, there remain \emph{significant challenges} in the use of stochastic optimization algorithms in practice.  The \emph{dirty secret} in the use of these algorithms is the tremendous computational costs required to \emph{tune} an algorithm for each application.  For large-scale, real-world systems, tuning an algorithm to solve a single problem might require weeks or months of effort on a supercomputer before the algorithm performs well.  To appreciate the consumption of energy to accomplish this, the authors of \cite{AsiDuch19} list multiple recent articles in which training a model for a single task requires thousands of CPU-days, and remark how each $10^4$ CPU-days is comparable to driving from Los Angeles to San Francisco with 50 Toyota Camrys.

One avenue for avoiding expensive tuning efforts is to employ \emph{adaptive} optimization algorithms.  Long the focus in the deterministic optimization community with widespread success in practice, such algorithms become significantly more difficult to design for the stochastic regime in which many modern problems reside, including those arising in large-scale ML and SP.

The purpose of this article is to summarize recent work and motivate continued research on the design and analysis of adaptive stochastic optimization methods.  In particular, we present an analytical framework---new for the context of adaptive deterministic optimization---that sets the stage for establishing convergence rate guarantees for adaptive stochastic optimization techniques.  With this framework in hand, we remark on important open questions related to how it can be extended further for the design of new methods.  We also discuss challenges and opportunities for their use in real-world systems.

\section{Background}

Many problems in ML and SP are formulated as optimization problems.  For example, given a data vector $y \in \R{m}$ from an unknown distribution, one often desires to have a vector of model parameters $x \in \R{n}$ such that a composite objective function $f : \R{n} \to \R{}$ is minimized, as in
\bequationNN
  \min_{x\in\R{n}}\ f(x),\ \text{where}\ f(x) := \E_y[\phi(x,y)] + \lambda \rho(x).
\eequationNN
Here, the function $\phi : \R{d+m} \to \R{}$ defines the \emph{data fitting} term $\E_y[\phi(x,y)]$.  For example, in supervised machine learning, the vector $y$ may represent the input and output from an unknown mapping and one aims to find $x$ to minimize the discrepancy between the output vector and the predicted value captured by $\phi$.  Alternatively, the vector $y$ may represent a noisy signal measurement and one may aim to find $x$ that filters out the noise to reveal the true signal.  The function $\rho : \R{n} \to \R{}$ with weight $\lambda \in [0,\infty)$ is included as a \emph{regularizer}.  This can be used to induce desirable properties of the vector $x$, such as sparsity, and/or to help avoid over-fitting a particular set of data vectors that are used when (approximately) minimizing $f$.  Supposing that instances of $y$ can be generated---one-by-one or in mini-batches, essentially \emph{ad infinitum}---the problem to minimize $f$ becomes a stochastic problem over $x$.

Traditional algorithms for minimizing $f$ are often very simple to understand and implement.  For example, given a solution estimate $x_k$, the well-known and celebrated stochastic gradient method (SG) \cite{RobbMonr51} computes the next estimate as $x_{k+1} \gets x_k - \alpha_k g_k$, where $g_k$ estimates the gradient of $f$ at $x_k$ by taking a random sample~$y_{i_k}$ and setting $g_k \gets \nabla_x \phi(x_k,y_{i_k}) + \lambda \nabla_x \rho(x_k)$ (or by taking a mini-batch of samples and setting $g_k$ as the average sampled gradient).  This value estimates the gradient since $\E_k[g_k] = \nabla f(x_k)$, where $\E_k[\cdot]$ conditions on the history of the behavior of the algorithm up to iteration $k \in \N{}$.  Under reasonable assumptions about the stochastic gradient estimates and with a \emph{prescribed} sequence of stepsize parameters $\{\alpha_k\}$, such an algorithm enjoys good convergence properties, which ensure that $\{x_k\}$ converges in probability to a minimizer, or at least a stationary point, of~$f$.

A practical issue in the use of SG is that the variance of the stochastic gradient estimates, i.e., $\E_k[\|g_k - \nabla f(x_k)\|_2^2]$, can be large, which inhibits the algorithm from attaining convergence rate guarantees on par with those for first-order algorithms in deterministic settings.  To address this, variance reduction techniques have been proposed and analyzed, such as those used in the algorithms SVRG, SAGA, and other methods \cite{CurtScheShi19,DefaBachLaco14,JohnZhan13,NguyLiuScheTaka17,schmidt2017minimizing}.  That said, SG and these variants of it are inherently \emph{nonadaptive} in the sense that each iteration involves a \emph{prescribed} number of data samples to compute $g_k$ in addition to a \emph{prescribed} sequence of stepsizes $\{\alpha_k\}$.  Determining which parameter values---defining the mini-batch sizes, stepsizes, and other factors---work well for a particular problem is a nontrivial task.  Tuning these parameters means that problems cannot be solved once; they need to be solved numerous times until reasonable parameter values are determined for future use on new data.

\section{Illustrative Example}

To illustrate the use of adaptivity in stochastic optimization, consider a problem of binary classification by logistic regression using the well-known MNIST dataset.  In particular, consider the minimization of a logistic loss plus an $\ell_2$-norm squared regularizer (with $\lambda = 10^{-4}$) in order to classify images as either of the digit 5 or not.

\begin{figure*}[ht]
  \bcenter
  \includegraphics[width=2in]{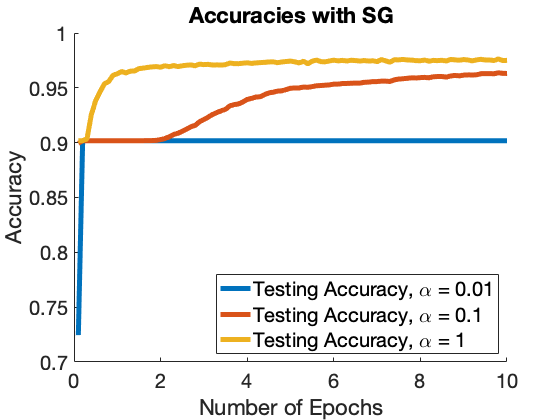}\qquad
  \includegraphics[width=2in]{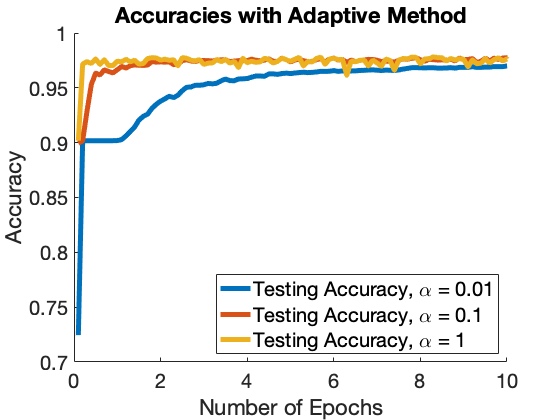}
  \ecenter
  \caption{SG vs.~an adaptive stochastic method on regularized logistic regression on MNIST.}
  \label{fig:mnist}
\end{figure*}

Employing SG with a mini-batch size of 64 and different fixed stepsizes, one obtains the plot of testing accuracy over 10 epochs given on the left in Figure~\ref{fig:mnist}.  One finds that for a stepsize of $\alpha_k = 1$ for all $k \in \N{}$, the model achieves testing accuracy around 98\%.  However, for a stepsize of $\alpha_k = 0.01$ for all $k \in \N{}$, the algorithm stagnates and never achieves accuracy much above 90\%.

By comparison, we also ran an adaptive method.  This method, like SG, begins with a mini-batch size of 64 and the stepsize parameter indicated in the plot on the right in Figure~\ref{fig:mnist}.  However, in each iteration, it checks the value of the objective (only over the current mini-batch) at the current iterate $x_k$ and trial iterate $x_k - \alpha_k g_k$.  If the mini-batch objective would not reduce sufficiently as a result of the trial step, then the step is not taken, the stepsize parameter is reduced by a factor, and the mini-batch size is increased by a factor.  This results in a more conservative stepsize with a more accurate gradient estimate in the subsequent iteration.  Otherwise, if the mini-batch objective would reduce sufficiently with the trial step, then the step is taken, the stepsize parameter is increased by a factor, and the mini-batch size is reduced by a factor.  Despite the data accesses required by this adaptive algorithm to evaluate mini-batch objective values in each iteration, the attained testing accuracy with all initializations competes with that attained by the best SG run.

This experiment demonstrates the potentially significant savings in computational costs offered by adaptive stochastic optimization methods.  While one might be able to achieve good practical performance with a nonadaptive method, it might only come after expensive tuning efforts.  By contrast, an adaptive algorithm can perform well without such expensive tuning.

\section{Framework for Analyzing Adaptive Deterministic Methods}\label{sec:general_determ}

Rigorous development of adaptive stochastic optimization methods requires a solid foundation in terms of convergence rate guarantees.  The types of adaptive methods that have enjoyed great success in the realm of deterministic optimization are (a) trust region, (b) line search, and (c) regularized Newton methods.  Extending these techniques to the stochastic regime is a highly nontrivial task.  After all, these methods traditionally require accurate function information at each iterate, which is what they use to adapt their behavior.  When an oracle can only return stochastic function estimates, comparing function values to make adaptive algorithmic decisions can be problematic.  In particular, when objective values are merely estimated, poor decisions can be made, and the combined effects of these poor decisions can be difficult to estimate and control.

As a first step toward showing how these challenges can be overcome, let us establish a general framework for convergence analysis for adaptive deterministic optimization.  This will lay a foundation for the framework that we present for adaptive stochastic optimization in \S\ref{sec:stoch-proc}.

The analytical framework presented here is new for the deterministic optimization literature.  A typical convergence analysis for adaptive deterministic optimization partitions the set of iterations into successful and unsuccessful ones.  Nonzero progress in reducing the objective function is made on successful iterations, whereas unsuccessful iterations merely result in an update of a model or algorithmic parameter to promote success in the subsequent iteration.  As such, a convergence rate guarantee results from a lower bound on the progress made in each successful iteration and a limit on the number of unsuccessful iterations that can occur between successful ones.  By contrast, in the framework presented here, the analysis is structured around a measure in which progress is made in \emph{all} iterations.

In the remainder of this section, we consider all three aforementioned types of adaptive algorithms under the assumption that $f$ is continuously differentiable with~$\nabla f$ being Lipschitz continuous with constant $L \in (0,\infty)$.  Each of these methods follows the general algorithmic framework that we state as Algorithm~\ref{alg:generic}.  The role of the sequence $\{\alpha_k\} \geq 0$ in the algorithm is to control the length of the trial steps $\{s_k(\alpha_k)\}$.  In particular, as seen in our discussion of each type of adaptive algorithm, one presumes that for a given model~$m_k$ the norm of $s_k(\alpha)$ is directly proportional to the magnitude of $\alpha$.  Another assumption---and the reason that we refer to it as a deterministic optimization framework---is that the models agree with the objective at least up to first-order derivatives, i.e., $m_k(x_k) = f(x_k)$ and $\nabla m_k(x_k) = \nabla f(x_k)$ for all $k \in \N{}$.

\begin{algorithm}[ht]
  \caption{Adaptive deterministic framework}
  \label{alg:generic}
  \begin{description}
    \item[Initialization]\ \\
    Choose constants $\eta \in (0,1)$, $\gamma \in (1,\infty)$, and $\overline\alpha \in (0,\infty)$.  Choose an initial iterate $x_0 \in \R{n}$ and stepsize parameter $\alpha_0 \in (0,\overline\alpha]$.
    \item[1. Determine model and compute step] \ \\
    Choose a local model $m_k$ of $f$ around $x_k$.  Compute a step $s_k(\alpha_k)$ such that the model reduction $m_k(x_k) - m_k(x_k+s_k(\alpha_k)) \geq 0$ is sufficiently large.
    \item[2. Check for sufficient reduction in $f$]\ \\
    Check if the reduction $f(x_k) - f(x_k+s_k(\alpha_k))$ is sufficiently large relative to the model reduction $m_k(x_k)-m_k(x_k+s_k(\alpha_k))$ using a condition parameterized by $\eta$.
    \item[3. Successful iteration]\ \\ 
    If sufficient reduction has been attained (along with other potential requirements), then set $x_{k+1} \gets x_k + s_k(\alpha_k)$ and $\alpha_{k+1} \gets 
\min\{\gamma\alpha_k,\overline\alpha\}$.
    \item[4. Unsuccessful iteration]\ \\  
    Otherwise, $x_{k+1} \gets x_k$ and $\alpha_{k+1} \gets \gamma^{-1} \alpha_k$.
    \item[5. Next iteration]\ \\
    Set $k \gets k + 1$.
  \end{description}
\end{algorithm}

Our analysis involves three central ingredients.  We define them here in such a way that they are easily generalized when we consider our adaptive stochastic framework later on.
\begin{itemize}
  \item $\{\Phi_k\} \geq 0$ is a sequence whose role is to measure progress of the algorithm.  The choice of this sequence may vary by type of algorithm and assumptions on $f$.
  \item $\{W_k\}$ is a sequence of indicators; specifically, for all $k \in \N{}$, if iteration $k$ is successful, then $W_k = 1$, and $W_k=-1$ otherwise. 
  \item $T_{\varepsilon}$, the \emph{stopping time}, is the index of the first iterate that satisfies a desired convergence criterion parameterized by ${\varepsilon}$. 
\end{itemize}
These quantities are not part of the algorithm itself, and therefore do not influence the iterates.  They are merely tools of the analysis.

At the heart of the analysis is the goal to show that the following condition holds.

\begin{condition}\label{cond.determ}
  The following statements hold with respect to $\{(\Phi_k,\alpha_k, W_k)\}$ and $T_{\varepsilon}$.
  \begin{enumerate}\setlength\itemsep{0em}
    \item There exists a scalar $\underline\alpha_\varepsilon \in (0,\infty)$ such that for each $k \in \N{}$ such that
    $\alpha_k \leq \gamma \underline\alpha_\varepsilon$, the iteration is guaranteed to be successful, i.e., $W_k=1$.  Therefore, $\alpha_k \geq \underline\alpha_\varepsilon$ for all $k \in \N{}$. 
    \item There exists a nondecreasing function $h_\varepsilon : [0,\infty) \to (0,\infty)$ and scalar $\Theta \in (0,\infty)$ such that, for all $k < T_{\epsilon}$,  $\Phi_k - \Phi_{k+1} \geq \Theta h_\varepsilon(\alpha_k)$.
      \end{enumerate}
\end{condition}

The goal to satisfy Condition~\ref{cond.determ} is motivated by the fact that, if it holds, it is trivial to derive (since $\Phi_k \geq 0$ for all $k \in \N{}$) that
\bequation\label{eq.Te}
  T_{\epsilon}\leq \frac{\Phi_0}{\Theta h_\varepsilon(\underline\alpha_\varepsilon)}.
\eequation
For generality, we have written $\underline\alpha_\varepsilon$ and $h_\varepsilon$ as parameterized by $\varepsilon$.  However, in the context of different algorithms, one or the other of these quantities may be independent of $\varepsilon$.  Throughout our analysis, we denote $f_* := \inf_{x\in\R{n}} f(x) > -\infty$.

\subsection{Classical trust region}\label{sec.det_tr}

In a classical trust region (TR) method, the model $m_k$ is chosen as at least a first-order accurate Taylor series approximation of $f$ at $x_k$, and the step $s_k(\alpha_k)$ is computed as a minimizer of $m_k$ in a ball of radius $\alpha_k$ centered at $x_k$.  In Step~2, the sufficient reduction condition is chosen as
\bequation\label{eq.tr_ratio}
  \frac{f(x_k)-f(x_k+s_k(\alpha_k))}{m(x_k)-m(x_k+s_k(\alpha_k))}\geq  \eta.
\eequation

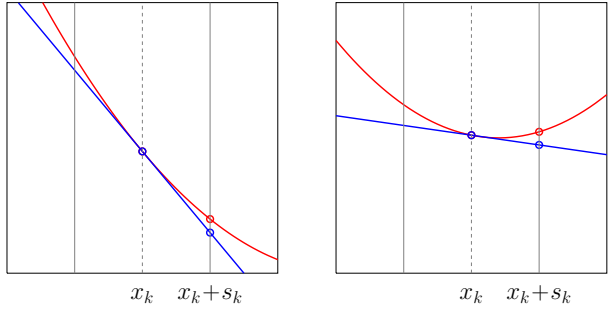
\begin{figure*}[ht]
  \bcenter
  \scalebox{0.45}{\input{det_good}}\qquad
  \scalebox{0.45}{\input{det_bad}}\ 
  \ecenter
  \caption{Illustration of successful (left) and unsuccessful (right) steps in a trust region method.}
  \label{fig:det_models}
\end{figure*}

Figure~\ref{fig:det_models} shows the need to distinguish between successful and unsuccessful iterations in a TR method.  Even though the model is (at least) first-order accurate, a large trust region radius may allow a large enough step such that the reduction predicted by the model does not well-represent the reduction in the function itself.  We contrast these illustrations later with situations in the stochastic setting, which are complicated by the fact that the model might be inaccurate regardless the size of the step.

For simplicity in our discussions, for iteration $k \in \N{}$ to be successful, we impose the additional condition that $\alpha_k\leq \tau \|\nabla f(x_k)\|$ for some suitably large constant $\tau \in (0,\infty)$.\footnote{Throughout the paper, consider all norms to be $\ell_2$.}  This condition is actually not necessary for the deterministic setting, but is needed in the stochastic setting to ensure that the trial step is not too large compared to the size of the gradient estimate.  We impose it now in the deterministic setting for consistency.

For this TR instance of Algorithm~\ref{alg:generic}, consider the first-order $\varepsilon$-stationarity stopping time
\bequation\label{eq.Te_tr}
  T_{\varepsilon} := \min\{k \in \N{} : \|\nabla f(x_k)\| \leq \varepsilon\},
\eequation
corresponding to which we define
\bequation\label{eq.Phi_tr}
  \Phi_k := \nu (f(x_k) - f_*) + (1-\nu)\alpha_k^2
\eequation
for some $\nu \in (0,1)$ (to be determined below).

Standard TR analysis involves two key results; see, e.g., \cite{NoceWrig06}.  First, while $\|\nabla f(x_k)\| > \varepsilon$,  if $\alpha_k$ is sufficiently small, then iteration~$k$ is successful, i.e., $W_k=1$; in particular, $\alpha_k \geq \underline\alpha_\varepsilon := c_1\varepsilon$ for all $k \in \N{}$ for some sufficiently small $c_1 \in (0,\infty)$ dependent on $L$, $\eta$, and $\gamma$.  In addition, if iteration~$k$ is successful, then the ratio condition~\eqref{eq.tr_ratio} and our imposed condition $\alpha_k\leq \tau \|\nabla f(x_k)\|$ yield
\bequationNN
  f(x_k) - f(x_{k+1}) \geq \eta c_2 \alpha_k^2
\eequationNN
for some $c_2 \in (0,\infty)$, meaning that
\bequationNN
  \Phi_k - \Phi_{k+1} \geq \nu \eta c_2 \alpha_k^2 - (1 - \nu)(\gamma^{2}-1)\alpha_k^2;
\eequationNN
otherwise, if iteration $k$ is unsuccessful, then
\bequationNN
  \Phi_k - \Phi_{k+1} = (1 - \nu)(1 - \gamma^{-2})\alpha_k^2.
\eequationNN
We aim to show in either case that $\Phi_k - \Phi_{k+1} \geq \Theta \alpha_k^2$, for some $\Theta>0$.  This can be done by choosing $\nu$ sufficiently close to $1$ such that
\bequationNN
 \nu \eta c_2 \geq (1 - \nu)(\gamma^{2} - \gamma^{-2}).
\eequationNN
In this manner, it follows from the observations above that Condition~\ref{cond.determ} holds with $h_\varepsilon(\alpha_k) := \alpha_k^2$ and $\Theta := (1-\nu)(1 - \gamma^{-2})$.  Thus, by \eqref{eq.Te}, the number of iterations required until a first-order $\varepsilon$-stationary point is reached satisfies
\bequationNN
  T_\varepsilon \leq \frac{\nu(f(x_0) - f_*) + (1 - \nu)\alpha^2_0}{(1-\nu)(1 - \gamma^{-2}) c_1^2\varepsilon^2}.
\eequationNN
This shows that $T_\varepsilon = \Ocal(\varepsilon^{-2})$.

\subsection{Classical line search}\label{sec.det_ls}

In a classical line search (LS) method, the model is again chosen as at least a first-order accurate Taylor series approximation of $f$ at $x_k$, with care taken to ensure it is convex so a minimizer of it exists. The trial step $s_k(\alpha_k)$ is defined as $\alpha_kd_k$ for some direction of sufficient descent $d_k$.  In Step~2, the sufficient reduction condition often includes the Armijo condition
\bequationNN
  f(x_k) - f(x_k + s_k(\alpha_k)) \geq -\eta \nabla f(x_k)^Ts_k(\alpha_k). 
\eequationNN
As is common, suppose $m_k$ is chosen and $d_k$ is computed such that, for a successful iteration, one finds for some $c_3 \in (0,\infty)$, dependent on $L$ and the angle between $d_k$ and $-\nabla f(x_k)$, that
\bequationNN
  f(x_k) - f(x_{k+1}) \geq \eta c_3 \alpha_k \|\nabla f(x_k)\|^2.
\eequationNN
Using common techniques, this can be ensured with $\alpha_k \geq \underline\alpha$ for all $k \in \N{}$ for some $\underline\alpha \in (0,\infty)$ dependent on $L$, $\eta$, and $\gamma$.  Similarly as in \S\ref{sec.det_tr}, let us also impose that $\|d_k\| \leq \beta \|\nabla f(x_k)\|$ for some suitably large $\beta \in (0,\infty)$.

For this LS instance of Algorithm~\ref{alg:generic}, for the stopping time $T_\varepsilon$ defined in \eqref{eq.Te_tr}, consider
\bequationNN
  \Phi_k := \nu (f(x_k) - f_*) + (1-\nu)\alpha_k\|\nabla f(x_k)\|^2. 
\eequationNN
If iteration $k$ is successful, then
\[
  \baligned
    & \Phi_k-\Phi_{k+1} \geq \nu \eta c_3 \alpha_k \|\nabla f(x_k)\|^2\\
    & -(1-\nu)( \gamma\alpha_k\|\nabla f(x_{k+1})\|^2-\alpha_k\|\nabla f(x_k)\|^2).
  \ealigned
\]
By Lipschitz continuity of $\nabla f$, it follows that
\[ 
\|\nabla f(x_{k+1})\| \leq (L\alpha_k\beta  + 1)\|\nabla f(x_k)\|.
\]
Squaring both sides and applying $\alpha_k\leq \bar \alpha$ yields
\[
 \|\nabla f(x_{k+1})\|^2  \leq  (L \bar \alpha\beta+1)^2 \|\nabla f(x_k)\|^2.
\]
Overall, if iteration $k$ is successful, then
\[
 \baligned
&\Phi_k-\Phi_{k+1} \geq \nu \eta c_3 \alpha_k \|\nabla f(x_k)\|^2 \\&-(1-\nu) 
 ((L \bar \alpha\beta+1)^2\gamma-1)\alpha_k \|\nabla f(x_k)\|^2.
\ealigned
\]
On unsuccessful iterations, one simply finds
\[
\Phi_k-\Phi_{k+1} \geq (1-\nu)(1-\gamma^{-1}) \alpha_k\|\nabla f(x_k)\|^2. 
\]
By selecting $\nu$ sufficiently close to $1$ such that
\[ \nu \eta c_3 \geq (1-\nu) 
 ((L\bar \alpha\beta+1)^2\gamma-\gamma^{-1}),
\]
one ensures, while $\|\nabla f(x_k)\| > \varepsilon$, that Condition~\ref{cond.determ} holds with $h_\varepsilon(\alpha_k) := \alpha_k\varepsilon^2 $ and $\Theta := (1-\nu)(1 - \gamma^{-1})$.  Thus, by \eqref{eq.Te}, the number of iterations required until a first-order $\varepsilon$-stationary point is reached satisfies
\bequationNN
  T_\varepsilon \leq \frac{\nu(f(x_0) - f_*) + (1 - \nu)\alpha_0\|\nabla f(x_0)\|^2}{(1-\nu)(1 - \gamma^{-1})\underline\alpha \varepsilon^2},
\eequationNN
which again shows that $T_\varepsilon = \Ocal(\varepsilon^{-2})$.

\subsection{Regularized Newton}\label{sec:rn_determ}

Regularized Newton methods have been popular in recent years due to their ability to offer optimal worst-case complexity guarantees for nonconvex smooth optimization over the class of practical second-order methods.  For example, in cubicly regularized Newton, the model is chosen as $m_k(x) = f(x_k) + \nabla f(x_k)^T(x-x_k) + \thalf(x-x_k)^T\nabla^2 f(x_k)(x-x_k) + \frac{1}{3\alpha_k} \|x-x_k\|^3$, and in Step~2 the imposed sufficient reduction condition has the same form as that in a trust region method, namely, the ratio condition~\eqref{eq.tr_ratio}.

As for the aforementioned line search method, one can show that if the Hessian of $f$ is Lipschitz continuous on a set containing the iterates, an iteration of cubicly regularized Newton will be successful if the stepsize parameter is sufficiently small; consequently, $\alpha_k \geq \underline\alpha$ for all $k \in \N{}$ for some $\underline\alpha \in (0,\infty)$ dependent on $L$, the Lipschitz constant for $\nabla^2 f$, $\eta$, and $\gamma$.  However, to prove the optimal complexity guarantee in this setting, one must account for the progress in a successful iteration as being dependent on the magnitude of the gradient at the \emph{next} iterate, not the current one.  (As we shall see, this leads to complications for the analysis in the stochastic regime.)  For this reason, and ignoring the trivial case when $\|\nabla f(x_0)\| \leq \varepsilon$, let us define
\bequationNN
  \baligned
    T_{\varepsilon} := \min\{k \in \N{}: \|\nabla f(x_{k+1})\| \leq \varepsilon\}
  \ealigned
\eequationNN
along with
\bequationNN
  \Phi_k := \nu (f(x_k) - f_*) + (1-\nu) \alpha_k \|\nabla f(x_k)\|^{3/2}.
\eequationNN
If iteration $k$ is successful, then
\bequationNN
  f(x_k)-f(x_{k+1}) \geq \eta c_4 \alpha_k \|\nabla f(x_{k+1})\|^{3/2}
\eequationNN
for some $c_4 \in (0,\infty)$ \cite{CartGoulToin11b}, meaning that
\bequationNN
  \baligned
    \Phi_k - \Phi_{k+1}
      \geq&\ \nu \eta c_4 \alpha_k \|\nabla f(x_{k+1})\|^{3/2} \\
        &\ + (1-\nu) \alpha_k \|\nabla f(x_k)\|^{3/2} \\
        &\ - (1-\nu) \alpha_{k+1} \|\nabla f(x_{k+1})\|^{3/2} \\
      \geq&\ \nu \eta c_4 \alpha_k \|\nabla f(x_{k+1})\|^{3/2} \\
        &\ - (1-\nu)\gamma \alpha_k \|\nabla f(x_{k+1})\|^{3/2};
  \ealigned
\eequationNN
otherwise, if iteration $k$ is unsuccessful, then
\bequationNN
  \Phi_k - \Phi_{k+1} \geq (1-\nu)(1-\gamma^{-1})\alpha_k \|\nabla f(x_{k+1})\|^{3/2}.
\eequationNN
Choosing $\nu$ sufficiently close to 1 such that
\bequationNN
  \nu \eta c_4 \geq (1 - \nu)(\gamma+ 1 - \gamma^{-1}),
\eequationNN
one ensures, while $\|\nabla f(x_{k+1})\| > \varepsilon$, that Condition~\ref{cond.determ} holds with $h_\varepsilon(\alpha_k) := \alpha_k \varepsilon^{3/2}$ and $\Theta := (1-\nu_\varepsilon)(1 - \gamma^{-1})$.  Thus, by \eqref{eq.Te}, the number of iterations required until a first-order $\varepsilon$-stationary point will be reached (in the next iteration) has
\bequationNN
  T_\varepsilon \leq \frac{\nu(f(x_0) - f_*) + (1 - \nu)\alpha_0\|\nabla f(x_0)\|^{3/2}}{(1-\nu)(1 - \gamma^{-1}) \underline\alpha \varepsilon^{3/2}},
\eequationNN
which shows that $T_\varepsilon = \Ocal(\varepsilon^{-3/2})$.

One can obtain the same result with a second-order TR method; see \cite{CurtRobiSama17}.  It should be said, however, that this method requires a more complicated mechanism for adjusting the stepsize parameter than that stated in Algorithm~\ref{alg:generic}.

\subsection{Additional examples}\label{sec:extras}

Our analysis so far has focused on the setting of having a stopping time based on a first-order stationarity criterion and no assumptions on $f$ besides first- and/or second-order Lipschitz continuous differentiability.  However, the framework can be extended for other situations as well.

For example, if one is interested in approximate second-order stationarity, then one can let
\bequationNN
  \baligned
    T_{\varepsilon} &:= \min\{k \in \N{} : \chi_k \leq \varepsilon\}, \\
    \text{where}\ \chi_k &:= \max\{\|\nabla f(x_k)\|, -\lambda_{\min}(\nabla ^2 f(x_k))\},
  \ealigned
\eequationNN
with $\lambda_{\min}(\cdot)$ denoting the minimum eigenvalue of its symmetric matrix argument.  One can show that if the model $m_k$ is chosen as a (regularized) second-order Taylor series approximation of~$f$ at~$x_k$, then for all of the aforementioned methods one obtains that $T_\varepsilon = \Ocal(\varepsilon^{-3})$. For TR, for example, one can derive 
this using 
\bequationNN
  \Phi_k := \nu (f(x_k) - f_*) + (1-\nu)\alpha_k^3. 
\eequationNN

On the other hand, if one is interested in analyzing specially the case of $f$ being convex, or even strongly convex, then one might consider
\bequationNN
  T_\varepsilon := \min\{k \in \N{} : f(x_k) - f_* \leq \varepsilon\}.
\eequationNN
In this case, improved complexity bounds can be obtained through other careful choices of $\{\Phi_k\}$.  For example, when
$f$ is convex and the LS method from \S\ref{sec.det_ls} is employed, one can let
\bequationNN
 \Phi_k := \nu \left (\frac{1}{\varepsilon}-\frac{1}{f(x_k) -f_*}\right )+ (1-\nu)\alpha_k. 
\eequationNN
Under the assumption that level sets of $f$ are bounded, one can show that $(f(x_{k+1}) - f_*)^{-1} - (f(x_k) - f_*)^{-1}$ is uniformly bounded below by a positive constant over all successful iterations, while for all $k \in \N{}$ one has $\alpha_k\geq \underline\alpha$.  Hence, for a suitable constant $\nu$, one can determine $h_\varepsilon(\alpha_k)$ and $\Theta$ to satisfy Condition~\ref{cond.determ}.  In this case, the function $h_\varepsilon$ does not depend on $\varepsilon$, but $\Phi_0 = \Ocal(\varepsilon^{-1})$, meaning that $T_\varepsilon=\Ocal(\varepsilon^{-1})$. 

Similarly, when $f$ is strongly convex and the LS method from \S\ref{sec.det_ls} is employed, consider
\bequationNN
\baligned
 \Phi_k& := \nu \left (\log\left(\frac{1}{\varepsilon}\right )-\log\left(\frac{1}{f(x_k) -f_*}\right )\right )\\ &+ (1-\nu)\log(\alpha_k). 
 \ealigned
\eequationNN
This time, $\log((f(x_{k+1}) - f_*)^{-1}) - \log((f(x_k) - f_*)^{-1})$ is uniformly bounded below by a positive constant over all successful iterations, and similar to the convex case one can determine $\nu$, $h_\varepsilon$, and $\Theta$ independent of $\varepsilon$ in order to show that $\Phi_0=\Ocal(\log(\varepsilon^{-1}))$ implies $ T_\varepsilon=\Ocal(\log(\varepsilon^{-1}))$. 

\section{Framework for Analyzing Adaptive Stochastic Methods}\label{sec:stoch-proc}

We now present a generalization of the framework introduced in the previous section that allows for the analysis of adaptive stochastic optimization methods.  This framework is based on the techniques proposed in \cite{cartis2018global}, which served to analyze the behavior of algorithms when derivative estimates are stochastic, but function values could be computed exactly.  It is also based on the subsequent work in \cite{BlanchetCartisMenickellyScheinberg2018}, which allows for function values to be stochastic as well.  For our purposes of providing intuition, we discuss a simplification of the framework, avoiding some technical details.  See \cite{cartis2018global,BlanchetCartisMenickellyScheinberg2018} for complete details.

As in the deterministic setting, let us define a generic algorithmic framework, which we state as Algorithm~\ref{alg:generic_stochastic}, that encapsulates multiple types of adaptive algorithms.  This algorithm has the same structure as Algorithm~\ref{alg:generic}, except it makes use of a stochastic model of $f$ to compute the trial step, and makes use of stochastic objective value estimates when determining whether sufficient reduction has been achieved.

\begin{algorithm}[ht]
  \caption{Adaptive stochastic framework}
  \label{alg:generic_stochastic}
  \begin{description}
    \item[Initialization]\ \\
    Choose $(\eta, \delta_1,\delta_2) \in (0,1)^3$, $\gamma\in (1, \infty)$, and $\overline\alpha \in (0,\infty)$.  Choose an initial iterate $x_0 \in \R{n}$ and stepsize parameter $\alpha_0 \in (0,\overline\alpha]$.
    \item[1. Determine model and compute step] \ \\
    Choose a stochastic model $m_k$ of $f$ around $x_k$, which satisfies some sufficient accuracy requirement with probability at least $1-\delta_1$.  Compute a step $s_k(\alpha_k)$ such that the model reduction $m_k(x_k) - m_k(x_k+s_k(\alpha_k)) \geq 0$ is sufficiently large.
    \item[2. Check for sufficient reduction]\ \\
    Compute estimates $\tilde f_k^0$  and $\tilde f_k^s$ of   $f(x_k)$ and  $f(x_k+s_k(\alpha_k))$, respectively, which satisfy some sufficient accuracy requirement with probability at least $1-\delta_2$.  Check if $\tilde f_k^0-\tilde f_k^s$ is sufficiently large relative to the model reduction $m_k(x_k)-m_k(x_k+s_k(\alpha_k))$ using a condition parameterized by~$\eta$.
    \item[3. Successful iteration]\ \\ 
    If sufficient reduction has been attained (along with other potential requirements), then set $x_{k+1} \gets x_k + s_k(\alpha_k)$ and $\alpha_{k+1} \gets 
\min\{\gamma\alpha_k,\overline\alpha\}$.
    \item[4. Unsuccessful iteration]\ \\  
    Otherwise, $x_{k+1} \gets x_k$ and $\alpha_{k+1} \gets \gamma^{-1} \alpha_k$.
    \item[5. Next iteration]\ \\
    Set $k \gets k + 1$.
  \end{description}
\end{algorithm}

Corresponding to Algorithm \ref{alg:generic_stochastic}, let $\{(\Phi_k, W_k)\}$ be a stochastic process such that $\{\Phi_k\} \geq 0$ for all $k \in \N{}$.  The sequences $\{\Phi_k\}$ and $\{W_k\}$ play similar roles as for our analysis of Algorithm~\ref{alg:generic}, but it should be noted that now each $W_k$ is a random indicator influenced by the iterate sequence $\{x_k\}$ and stepsize parameter sequence $\{\alpha_k\}$, which are themselves  stochastic processes.

Let $\Fcal_k$ denote the $\sigma$-algebra generated by all stochastic processes within Algorithm~\ref{alg:generic} at the {\em beginning} of iteration $k$.  Roughly, $\Fcal_k$ it is generated by $\{(\Phi_j,\alpha_j, x_j)\}_{j=0}^k$.  Note that this includes $\{(\Phi_j,\alpha_j, x_j)\}_{j=0}^{k-1}$, but does not include $W_k$ as this random variable depends on what happens {\em during} iteration $k$. We then let $T_\varepsilon$ denote a family of stopping times for $\{(\Phi_k,W_k)\}$ with respect to $\{\Fcal_k\}$ parameterized by $\varepsilon \in (0,\infty)$.

The goal of our analytical framework is to derive an upper bound for the expected stopping time $\E[T_{\varepsilon}]$ under various assumptions on the behavior of $\{(\Phi_k,W_k)\}$ and on the objective~$f$.  At the heart of the analysis is the goal to show that the following condition holds, which can be seen as a generalization of Condition~\ref{cond.determ}.

\begin{condition}\label{cond.stoch}
  The following statements hold with respect to $\{(\Phi_k,\alpha_k, W_k)\}$ and $T_{\varepsilon}$.
  \begin{enumerate}\setlength\itemsep{0em}

    \item There exists a scalar $\underline\alpha_\varepsilon \in (0,\infty)$ such that, conditioned on the event that $\alpha_k \leq \underline\alpha_\varepsilon$, one has $W_k = 1$ with probability $1 - \delta > \thalf$, conditioned on $\Fcal_k$.   
    $($This means that, if it becomes sufficiently small, one is more likely to see an increase in the stepsize parameter than a decrease in it.$)$
    \item There exists a nondecreasing function $h_\varepsilon : [0,\infty) \to (0,\infty)$ and scalar $\Theta \in (0,\infty)$ such that, for all $k < T_{\varepsilon}$, the conditional expectation of $\Phi_k - \Phi_{k+1}$ with respect to $\Fcal_k$ is at least $\Theta h_\varepsilon(\alpha_k)$; specifically,
  \end{enumerate}
  \bequationNN
    \mathds{1}_{\{k < T_{\varepsilon}\}} \E[\Phi_{k+1} | \Fcal_k] \leq \mathds{1}_{\{k < T_{\varepsilon}\}} (\Phi_k - \Theta h_\varepsilon(\alpha_k)).
  \eequationNN
\end{condition}

Whereas Condition~\ref{cond.determ} requires that the stepsize parameter $\alpha_k$ remains above a lower bound and the reduction $\Phi_k - \Phi_{k+1}$ is nonnegative with certainty, Condition~\ref{cond.stoch} allows more flexibility.  In particular, it says that until the stopping time is reached, one \emph{tends} to find $\Phi_k - \Phi_{k+1}$ sufficiently large, namely, at least $\Theta h_\varepsilon(\alpha_k)$ for each $k < T_\varepsilon$.  In addition, the stepsize parameter is allowed to fall below the threshold~$\underline\alpha_\varepsilon$; in fact, it can become arbitrarily small.  That said, if $\alpha_k\leq \underline\alpha_\varepsilon$, then one \emph{tends} to find $\alpha_{k+1} > \alpha_k$.

With this added flexibility, one can still  prove complexity guarantees.  Intuitively, the reduction $\Phi_k - \Phi_{k+1}$ is at least the deterministic amount $\Theta h_\varepsilon(\underline\alpha_\varepsilon)$ often enough that one is able to bound the total number of such occurrences (since $\{\Phi_k\} \geq 0$).  The following theorem (see Theorem~2.2 in \cite{BlanchetCartisMenickellyScheinberg2018}) bounds the expected stopping time in terms of a deterministic value.

\begin{thm}\label{thm:renewal_reward_stop_time}
  If Condition~\ref{cond.stoch} holds, then
  \bequationNN
    \E[T_{\varepsilon}] \leq \frac{1-\delta}{1-2\delta} \cdot
    \frac{\Phi_0}{\Theta h_\varepsilon(\underline\alpha_\varepsilon)} 
    + 1.
  \eequationNN
\end{thm}

In \S\ref{sec.trust_region}--\ref{sec.line_search}, we summarize how this framework has been applied to analyze the behavior of stochastic TR and LS methods.  In each case, the keys to applying the framework are determining how to design the process $\{\Phi_k\}$ as well as specify details of  Algorithm \ref{alg:generic_stochastic}  to ensure that Condition~\ref{cond.stoch} holds.  For these aspects, we first need to describe different adaptive accuracy requirements for stochastic function and derivative estimates that might be imposed in Steps~1 and 2, as well the techniques that have been developed to ensure these requirements.

\subsection{Error Bounds for Stochastic Function and Derivative Estimates}\label{sec:bounds}

In this section, we describe various types of conditions that one may require in an adaptive stochastic optimization algorithm when computing objective function, gradient, and Hessian estimates.  These conditions have been used in some previously proposed adaptive stochastic optimization algorithms \cite{BlanchetCartisMenickellyScheinberg2018, RByrd_etal_2012, cartis2018global, PaquetteScheinberg2018}.

Let us remark in passing that one does not necessarily need to employ sophisticated error bound conditions in stochastic optimization to achieve improved convergence rates.  For example, 
in the case of minimizing strongly convex~$f$, the linear rate of convergence of \emph{gradient descent} can be emulated by an SG-like method if the mini-batch size grows exponentially \cite{Schmidt,  PasuGlynGhosFate18}.  However, attaining similar improvements in the (not strongly) convex and nonconvex settings has proved elusive.  Moreover,  while the stochastic estimates improve with the progress of such an algorithm, this improvement is based on \emph{prescribed} parameters, and hence the algorithm is not adaptive in our desirable sense.  Hence, one still needs to tune such algorithms for each application.

Returning to our setting, in the types of error bounds presented below, for a given $f : \R{n} \to \R{}$ and $x \in \R{n}$, let $\ftilde(x) $, $g(x)$, and $H(x)$ denote stochastic approximations of $f(x)$, $\nabla f(x)$, and $\nabla^2 f(x)$, respectively.

\begin{itemize}\setlength\itemsep{0em}
  \item \textbf{Taylor-like Conditions}.  Corresponding to a norm $\|\cdot\|$, let $\Bmbb(x_k,\Delta_k)$ denote a ball of radius~$\Delta_k$ centered at $x_k$.  If the function, gradient, and Hessian estimates satisfy
  \bsubequations\label{eq:fully-linear}
    \begin{align}
      |\ftilde(x_k) - f(x_k)| &\leq \kappa_f \Delta_k^2, \label{eq:fully-linear_f} \\
      \|g(x_k) - \nabla f(x_k)\| &\leq \kappa_g \Delta_k, \label{eq:fully-linear_g} \\ \text{and}\ \ 
      \|H(x_k) - \nabla^2 f(x_k)\| &\leq \kappa_H \label{eq:fully-linear_H}
    \end{align}
  \esubequations
  for some nonnegative scalars $(\kappa_f,\kappa_g,\kappa_H)$, then the model $m_k(x) = \ftilde(x_k) + g(x_k)^T(x - x_k) + \thalf (x - x_k)^TH(x_k)(x - x_k)$ gives an approximation of $f$ within $\Bmbb(x_k,\Delta_k)$ that is comparable to that given by an accurate first-order Taylor series approximation (with error dependent on $\Delta_k$).  Similarly, if~\eqref{eq:fully-linear} holds with the right-hand side values replaced by $\kappa_f \Delta_k^3$, $\kappa_g \Delta_k^2$, and $\kappa_H\Delta_k$, respectively, then $m_k$ gives an approximation of~$f$ that is comparable to that given by an accurate second-order Taylor series approximation.  In a stochastic setting when unbiased estimates of $f(x_k)$, $\nabla f(x_k)$, and $\nabla^2 f(x_k)$ 
  can be computed, such conditions can be ensured, with some sufficiently high probability $1-\delta$, by sample average approximations using a sufficiently large number of samples. For example, to satisfy  \eqref{eq:fully-linear_g}  with probability $1-\delta$, the sample size for computing $g(x_k)$ can be $\Omega \(\tfrac{V_g}{\kappa_g^2\Delta_k^2}\)$, where $V_g$ is the variance of the stochastic gradient estimators corresponding to a mini-batch size of~1.  Here, $\Omega$-notation hides dependence on $\delta$, which is weak if $\|g(x_k)-\nabla f(x_k)\|$ is bounded.
    
  \item \textbf{Gradient Norm Condition.}  If, for some $\theta \in [0,1)$, one has that
  \bequation\label{eq:norm-cond}
    \|g(x_k) - \nabla f(x_k)\| \leq \theta \|\nabla f(x_k)\|,
  \eequation
  then we say that $g(x_k)$ satisfies the gradient norm condition at $x_k$.  Unfortunately, verifying the gradient norm condition at $x_k$ requires knowledge of $\|\nabla f(x_k)\|$, which makes it an impractical condition.  In \cite{RByrd_etal_2012}, a heuristic is proposed that attempts to approximate the sample size for which the gradient norm condition holds.  More recently, in \cite{Bollapragada2017}, the authors improve upon the gradient norm condition by introducing an \emph{angle condition}, which in principle allows smaller sample set sizes to be employed.  However, again, the angle condition requires a bound in terms of $\| \nabla f(x_k)\|$, for which a heuristic estimate needs to be employed.
  
  Instead of employing the unknown quantity $\|\nabla f(x_k)\|$ on the right hand side of the norm and angle conditions, one can substitute $\varepsilon \in (0,\infty)$, the desired stationarity tolerance.  In this manner, while $\|\nabla f(x_k)\| > \varepsilon$ and $\|g_k - \nabla f(x_k)\| \leq \theta \varepsilon$, the norm condition is satisfied.  This general idea has been exploited recently in several articles (e.g., \cite{Tripuraneni2017,Mahoney}), where gradient and Hessian estimates are computed based on large-enough numbers of samples, then assumed to be accurate in \emph{every} iteration (with high probability) until an $\epsilon$-stationary solution is reached.  While strong iteration complexity guarantees can be proved for such algorithms, these approaches are too conservative to be competitive with truly stochastic algorithms.

  \item \textbf{Stochastic Gradient Norm Conditions.} 
  Consider again the conditions in~\eqref{eq:fully-linear}, although now consider the specific setting of defining $\Delta_k := \alpha_k\|g(x_k)\|$ for all $k \in \N{}$.  This condition, which involves a bound comparable to \eqref{eq:norm-cond} when $g(x_k) = \nabla f(x_k)$, is particularly useful in the context of LS methods.  It is possible to impose it in probability since, unlike \eqref{eq:norm-cond}, it does not require knowledge of $\|\nabla f(x_k)\|$.  In this case,  to satisfy  \eqref{eq:fully-linear} for $\Delta_k := \alpha_k\|g(x_k)\|$ with probability $1-\delta$, the sample size for computing $g(x_k)$ need only be $\Omega \(\tfrac{V_g}{\kappa_g^2\alpha_k^2 \|g(x_k)\|^2}\)$.   While $g(x_k)$ is not known when the sample size for computing it is chosen, one can use a simple loop that guesses the value of $\|g(x_k)\|$, then iteratively increases the number of samples as needed; see \cite{cartis2018global,PaquetteScheinberg2018}. 
\end{itemize}
 
Note that all of the above conditions can be adaptive in terms of the progress of an algorithm, assuming that $\{\Delta_k\}$ and/or $\{\|g(x_k)\|\}$ vanish as $k \to \infty$.  However, this behavior does not have to be monotonic, which is a benefit of adaptive stepsizes and accuracy requirements.


\subsection{Stochastic trust region}\label{sec.trust_region}

The idea of designing stochastic TR methods has been attractive for years, even before the recent explosion of efforts on algorithms for stochastic optimization; see, e.g., \cite{ChangHongWan}.  This is due to the impressive practical performance that TR methods offer, especially when (approximate) second-order information is available.  Since TR methods typically compute trial steps by minimizing a quadratic model of the objective in a neighborhood around the current iterate, they are inherently equipped to avoid certain spurious stationary points that are not local minimizers.  In addition, by normalizing the step length by a trust region radius, the behavior of the algorithm is kept relatively stable.  Indeed, this feature allows TR methods to offer stability even in the non-adaptive stochastic regime; see \cite{CurtScheShi19}.

\begin{figure*}[ht]
  \bcenter
  \scalebox{0.45}{\input{good_good}}\ 
  \scalebox{0.45}{\input{good_bad}}\ 
  \scalebox{0.45}{\input{bad_good}}\ 
  \scalebox{0.45}{\input{bad_bad}}
  \ecenter
  \caption{Illustration of ``good'' and ``bad'' models and estimates in a stochastic trust region method: (first) good model and good estimates; (second) good model and bad estimates; (third) bad model and good estimates; (fourth) bad model and bad estimates.}
  \label{fig:models}
\end{figure*}
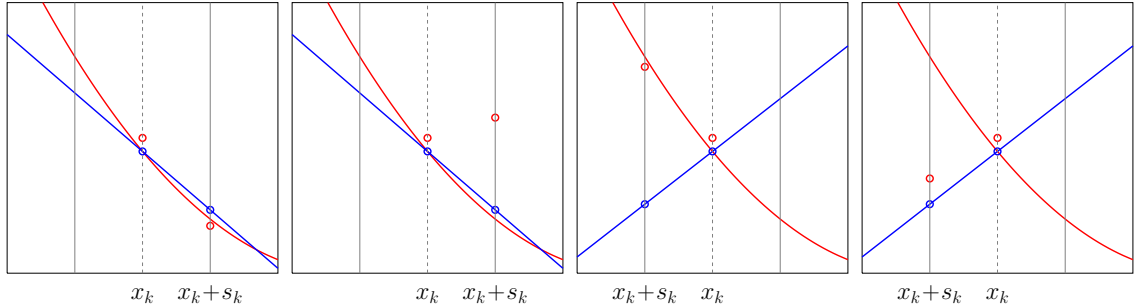


However, it has not been until the last couple years that researchers have been able to design stochastic TR methods that can offer strong expected complexity bounds, which are essential for ensuring that the practical performance of such methods can be competitive over broad classes of problems.  The recently proposed algorithm known as STORM, introduced in \cite{chen2018stochastic}, achieves such complexity guarantees by requiring the Taylor-like conditions \eqref{eq:fully-linear} to hold with probability at least $1-\delta$, conditioned on $\Fcal_k$.\footnote{In what follows, all accuracy conditions are assumed to hold with some probability, conditioned on $\Fcal_k$.  We omit mention of this conditioning for brevity.}  In particular, in Step 1 of Algorithm \ref{alg:generic_stochastic} a model $m_k(x) := \tilde f(x_k) + g(x_k)^T (x - x_k) + \tfrac12 (x - x_k)^T H(x_k) (x - x_k)$ is computed with components satisfying \eqref{eq:fully-linear} with probability $1-\delta_1$ (where $\Delta_k := \alpha_k$ for all $k \in \N{}$), then the step $s_k(\alpha_k)$ is computed by minimizing $m_k$ (approximately) within a ball of radius $\alpha_k$. In Step 2, the estimates 
$\tilde f_k^0$ and $\tilde f_k^s$ are computed to  satisfy \eqref{eq:fully-linear_f} with probability $1-\delta_2$. The imposed sufficient reduction condition is 
\bequationNN
 \dfrac{\tilde f_k^0-\tilde f_k^s}{m_k(x_k)-m_k(x_k+s_k)}\geq \eta. 
\eequationNN
An iteration is successful if the above holds and $ \| g_k\| \geq \tau \alpha_k$ for some user-defined $\tau \in (0,\infty)$. 

For brevity, we omit some details of the algorithm.  For example, the constant $\kappa_f$ in \eqref{eq:fully-linear} cannot be too large, whereas $\kappa_g$ and $\kappa_H$ can be arbitrarily large.  Naturally, the magnitudes of these constants affects the constants in the convergence rate; see~\cite{BlanchetCartisMenickellyScheinberg2018} for further details.

Based on the stochastic process generated by Algorithm \ref{alg:generic_stochastic}, let us define, as in the deterministic setting, $T_\varepsilon$ and $\{\Phi_k\}$ by \eqref{eq.Te_tr} and \eqref{eq.Phi_tr}, respectively.  It is shown in \cite{BlanchetCartisMenickellyScheinberg2018} that for sufficiently small constants $\delta_1$, $\delta_2$, and $\Theta$ (independent of~$\varepsilon$), Condition~\ref{cond.stoch} holds with 
\begin{equation*}
\underline {\alpha}_\varepsilon := \zeta\varepsilon \ \ \text{and}\ \ h(\alpha_k) := \alpha_k^2
\end{equation*}
for some positive constants $\zeta$ and $\Theta$ that depend on the algorithm parameters and properties of $f$, but not on $\varepsilon$.  The probability $1-\delta$ that arises in Condition~\ref{cond.stoch} is at least $(1-\delta_1)(1-\delta_2)$. Thus, by
Theorem \ref{thm:renewal_reward_stop_time}, the expected complexity of reaching a first-order $\varepsilon$-stationary point is at most $\Ocal(\varepsilon^{-2})$, which matches the complexity of the deterministic version of the algorithm up to the factor dependent on $(1-\delta_1)(1-\delta_2)$.

Let us give some intuition of how Condition~\ref{cond.stoch} is ensured. Let us say that the model $m_k$ is ``good'' if its components satisfy \eqref{eq:fully-linear} and ``bad'' otherwise. Similarly, the estimates $f_k^0$ and $f_k^s$ are ``good'' if they satisfy  \eqref{eq:fully-linear_f} and ``bad'' otherwise.  By construction of Steps~1 and 2 of the algorithm, $m_k$ is ``good'' with probability at least $1-\delta_1$ and $f_k^0$ and $f_k^s$ are ``good'' with probability at least $1-\delta_2$.  Figure~\ref{fig:models} illustrates the four possible outcomes.  When both the model and the estimates are good, the algorithm essentially behaves as its deterministic counterpart; in particular, if $\alpha_k\leq \underline \alpha_\varepsilon$ and $k < T_\varepsilon$, then the $k$-th iteration is successful and the reduction $\Phi_k - \Phi_{k+1}$ is sufficiently large.  If the model is bad and the estimates are good, or if the model is good and the estimates are bad, then the worst case (depicted in Figure \ref{fig:models}) is that the step is deemed unsuccessful, even though $\alpha_k\leq \underline \alpha_\varepsilon$.  This shows that the stepsize parameter can continue to decrease, even if it is already small.  Finally, if both the model and estimates are bad, which happens with probability at most $\delta_1\delta_2$, then it is possible that the iteration will be deemed successful despite the fact that $\Phi_{k+1} > \Phi_k$.  (Recall that this cannot occur in the deterministic setting, where $\{\Phi_k\}$ decreases monotonically.)  The key step in showing that $\Phi_{k+1} < \Phi_k$ \emph{in expectation} is to establish that, in each iteration, the possible decrease of this measure is proportional to any possible increase, and thus by ensuring that $\delta_1\delta_2$ is sufficiently small and $(1-\delta_1)(1-\delta_2)$ is sufficiently large, one can ensure a desired reduction in expectation. 

The same TR algorithm can be employed with minor modifications to obtain good expected complexity properties with respect to achieving second-order $\varepsilon$-stationarity.  In this case, the requirements on the estimates need to be stronger; in particular, \eqref{eq:fully-linear} has to be imposed with right-hand side values $\kappa_f \alpha_k^3$, $\kappa_g \alpha_k^2$, and $\kappa_H\alpha_k$, respectively.  In this case, Condition~\ref{cond.stoch}  holds for $h(\alpha_k)= \alpha_k^3$.  Hence, the expected complexity of reaching an second-order $\varepsilon$-stationary point by Algorithm \ref{alg:generic_stochastic} is bounded by ${\cal O}(\varepsilon^{-3})$, which similarly matches the deterministic complexity. 

\subsection{Line Search Methods}\label{sec.line_search}

A major disadvantage of an algorithm such as SG is that one is very limited in the choice of stepsize sequence that can be employed to adhere to the theoretical guidelines.  One would like to be able to employ a type of line search, as has been standard practice throughout the history of research on deterministic optimization algorithms.  However, devising line search methods for the stochastic regime turns out to be extremely difficult.  This is partially due to the fact that, unlike TR methods, LS algorithms employ steps that are highly influenced by the norm of the gradient $\|\nabla f(x_k)\|$, or in the stochastic regime influenced by $\|g(x_k)\|$.  Since the norm of the gradient estimate can vary dramatically from one iteration to the next, the algorithm needs to have a carefully controlled stepsize selection mechanism to ensure convergence.
 
In  \cite{Bollapragada2017} and \cite{Schmidt}, two backtracking line search methods were proposed that use different heuristic sample size strategies when computing gradient and function estimates. In both cases, the backtracking is based on the  Armijo condition applied to function estimates that are computed on the same batch as the gradient estimates. A different type of LS method that uses a probabilistic Wolfe condition for choosing the stepsize was proposed in \cite{Prob_line_search}, although this approach possesses no known theoretical guarantees.

In~\cite{Vasanietal2019}, the authors argue that with overparametrization of deep neural networks (DNNs), the variance of the stochastic gradients tends to zero near stationarity points. Under this assumption, the authors analyzed a stochastic LS method, which has good practical performance in some cases.  However, not only is such an assumption unreasonably strong, even for DNNs, but it also does not extend to the fully stochastic setting, since in that setting the assumption would imply zero generalization error at the solution---an ideal situation, but hardly realistic in general.

Here we summarize results in \cite{PaquetteScheinberg2018}, where an LS method with an adaptive sample size selection mechanism is proposed and  complexity bounds are provided.  This method can be described as a particular case of Algorithm \ref{alg:generic_stochastic}. As in the deterministic case, a stochastic model $m_k$ is chosen in Step~1 and $s(\alpha_k)=-\alpha_k d_k$, where $d_k$ makes an obtuse angle with the gradient estimate $g(x_k)$.  The sufficient reduction in Step~2 is based on the estimated Armijo condition
\bequationNN
  \ftilde_k^0 - \ftilde_k^s \geq -\eta g(x_k)^Ts_k(\alpha_k). 
\eequationNN
The algorithm requires that the components of the model $m_k$ satisfy \eqref{eq:fully-linear} with probability at least $1-\delta_1$, and that the estimates $\tilde f_k^0$ and $\tilde f_k^s$ satisfy \eqref{eq:fully-linear_f} with probability at least $1-\delta_2$.  Here, it is critical that \eqref{eq:fully-linear_f} not be imposed with the right-hand side being $\alpha_k\|g_k\|^2$ (even though the deterministic case might suggest this as being appropriate) since this quantity can vary uncontrollably from one iteration to the next.  To avoid this issue, the approach defines an additional control sequence $\{\Delta_k\}$ used for controlling the accuracy of $\tilde f_k^0$ and $\tilde f_k^s$.   Intuitively, for all $k \in \N{}$, the value $\Delta_k^2$ is meant to approximate $\alpha_k\|\nabla f(x_k)\|^2$, which, as seen in the deterministic case, is the desired reduction in $f$ if iteration $k$ is successful.  This control sequence needs to be set carefully.   The first value in the sequence is set arbitrarily, with subsequent values set as follows.  If iteration~$k$ is unsuccessful, then $\Delta_{k+1} \gets \Delta_k$.  Otherwise, if iteration~$k$ is successful, then one checks whether the step is \emph{reliable} in the sense that
  the accuracy parameter is sufficiently small, i.e.,
  \bequation\label{eq.reliable}
    \alpha_k \|g(x_k)\|^2 \geq \Delta_k^2.
  \eequation
  If \eqref{eq.reliable} holds, then one sets $\Delta_{k+1} \gets \sqrt{\gamma} \Delta_k$; otherwise, one sets $\Delta_{k+1} \gets \sqrt{\gamma^{-1}} \Delta_k$ to promote reliability of the step in the subsequent iteration.
     Using $\{\Delta_k\}$ defined in this manner, an additional bound is imposed on the variance of 
     the objective value estimates. For all $k \in \N{}$, one requires
\bequationNN
\begin{aligned}
 \max&\{ \E| \tilde f_k^0  - f(x_k)|^2, \E| \tilde f_k^s  - f(x_k+s_k(\alpha_k))|^2\}\\
  & \leq \max\{\kappa_f \alpha_k^2 \|\nabla f(x_k)\|^2, \kappa_f \Delta_k^4\}.\\
 \end{aligned}
\eequationNN
Note that because of the $\max$ in the right hand side of this inequality, and because $\Delta_k$ is a known value, it is not necessary to know $\|\nabla f(x_k)\|^2$ in order to impose this condition. Also note that this condition is stronger than imposing  the Taylor-like condition \eqref{eq:fully-linear_f} with some probability less than $1$, because the bound on expectation does not allow $| \tilde f_k^0  - f(x_k)|$ to be arbitrarily large with positive probability, while \eqref{eq:fully-linear_f} allows it, and thus is more tolerant to outliers.

For analyzing this LS instance of Algorithm \ref{alg:generic_stochastic}, again let $T_{\varepsilon}$ be as in \eqref{eq.Te_tr}.  However, now let 
 \bequationNN
   \Phi_k := \nu (f(x_k)-f_*)+ (1-\nu) (\tfrac{\alpha_k}{L^2}\|\nabla    f(x_k)\|^2+   \eta\Delta_k^2).
 \eequationNN
Using a strategy similar to that for STORM combined with  the logic of \S\ref{sec.det_ls} it is shown in \cite{PaquetteScheinberg2018} that Condition~\ref{cond.stoch} holds with
\[
  h(\alpha_k) = \alpha_k\varepsilon^2. 
\]



The expected complexity of this stochastic LS method has been analyzed for minimizing convex and strongly convex $f$ using modified definitions for $\Phi_k$ and $T_\varepsilon$ as described in \S\ref{sec:extras}; see   \cite{PaquetteScheinberg2018}.


\subsection{Stochastic regularized Newton}\label{sec:prior_arc}

As we have seen, stochastic TR and LS algorithms have been developed that fit into the adaptive stochastic framework that we have described, showing that they can achieve expected complexity guarantees on par with their deterministic counterparts.  However, \emph{neither of these types of algorithms achieves complexity guarantees that are optimal in the deterministic regime}.  

Cublicly regularized Newton \cite{CartGoulToin11a,CartGoulToin11b,NesterovPolyak2006}, described and analyzed in \S\ref{sec:rn_determ}, enjoys optimal convergence rates for second-order methods for minimizing nonconvex functions.  We have shown how our adaptive deterministic framework gives  a  $\Ocal(\varepsilon^{-3/2})$ complexity bound for this method  for achieving first-order $\varepsilon$-stationarity, in particular, to achieve $\|\nabla f(x_{k+1})\| \leq \varepsilon$.  There has also been some  work that proposes stochastic and randomized versions of cubic regularization methods \cite{Tripuraneni2017,Mahoney}, but these impose strong conditions on the accuracy of the function, gradient, and Hessian estimates that are equivalent to using $\varepsilon$ (i.e., the desired accuracy threshold) in place of $\Delta_k$ for all $k \in \N{}$ in~\eqref{eq:fully-linear}.  Thus, these approaches essentially reduce to sample average approximation with very tight accuracy tolerances and are not adaptive enough to be competitive with truly stochastic algorithms in practice.  For example, in \cite{Tripuraneni2017}, no adaptive stepsizes or batch sizes are employed, and in \cite{Mahoney} only Hessian approximations are assumed to be stochastic.  Moreover, the convergence analysis is performed under the assumption that the estimates are sufficiently accurate (for a given $\varepsilon$) in {\em every} iteration.  Thus, essentially, the analysis is reduced to that in the deterministic setting and applies only as long as no iteration fails to satisfy the accuracy condition. Hence,  a critical open question is whether one can extend the framework described here (or another approach) to develop and  analyze a stochastic algorithm  that achieves \emph{optimal} deterministic complexity. 

The key difficulty in extending the analysis described in \S\ref{sec:rn_determ} to the stochastic regime is the definition of the stopping time. For Theorem \ref{thm:renewal_reward_stop_time} to hold, $T_\varepsilon$ has to be a stopping time   with respect to $\{\Fcal_k\}$. However, with $s_k(\alpha_k)$ being random, $x_{k+1}$ is not measurable in $\{\Fcal_k\}$; hence, $T_\varepsilon$, as it is defined in the deterministic setting, is not a valid stopping time in the stochastic regime.  A different definition is needed that would be agreeable with the analysis in the stochastic setting. 

Another algorithmic framework that enjoys optimal complexity guarantees is the Trust Region Algorithm with Contractions and Expansions (TRACE) \cite{CurtRobiSama17}.  This algorithm borrows much from the traditional TR methodology also followed by STORM, but incorporates a few algorithmic variations that reduces the complexity from $\Ocal(\varepsilon^{-2})$ to $\Ocal(\varepsilon^{-3/2})$ for achieving first-order $\varepsilon$-stationarity.  It remains an open question whether one can employ the adaptive stochastic framework to analyze a stochastic variant of TRACE, the main challenge being that TRACE involves a relatively complicated strategy for updating the stepsize parameter.  Specifically, deterministic TRACE requires knowledge of the exact Lagrange multiplier of the trust region constraint at a solution of the step computation subproblem.  If the model $m_k$ is stochastic and the subproblem is solved only approximately, then it remains open how to maintain the optimal (expected) complexity guarantee. In addition, the issue of determining the correct stopping time in the stochastic regime is also open for this method. 

\section{Other possible extensions}

We have shown that the analytical framework for analyzing adaptive stochastic optimization algorithms presented in \S\ref{sec:stoch-proc} has offered a solid foundation upon which stochastic TR and stochastic LS algorithms have been proposed and analyzed.  We have also shown the challenges and opportunities for extending the use of this framework for analyzing algorithms whose deterministic counterparts have optimal complexity.

Other interesting questions remain to be answered.  For example, great opportunities exist for the design of error bound conditions other than those mentioned in \S\ref{sec:bounds}, especially when it comes to bounds that are tailored for particular problem settings.  While improved error bounds might not lead to improvements in iteration complexity, they can have great effects on the work complexity of various algorithms, which translates directly into performance gains in practice.

The proposed analytical framework might also benefit from extensions in terms of the employed algorithmic parameters.  For example, rather than using a single constant $\gamma$ when updating the stepsize parameter, one might consider different values for increases vs.~decreases, and when different types of steps are computed.  This will allow for improved bounds on the accuracy probability tolerances $\delta_1$ and $\delta_2$.

Finally, numerous open questions remain in terms of how best to implement adaptive stochastic algorithms in practice.  For different algorithm instances, practitioners need to explore how best to adjust mini-batch sizes and stepsizes so that one can truly achieve good performance without wasteful tuning efforts.  One hopes that with additional theoretical advances, these practical questions will become easier to answer.

\bibliographystyle{plain}
\bibliography{fec,Katya,references,white_paper}

\ifappendix

\onecolumn

\section*{Appendix}

\bitemize
  \item Analysis in \S\ref{sec.det_tr}.  The lower bound on $\{\alpha_k\}$ can be seen in the proof of Theorem~4.5 in Nocedal \& Wright (2006).  In particular, under the assumption that $\|\nabla f(x_k)\| > \varepsilon$ for all $k \in \N{}$, one can derive a lower bound on the trust region radius of the form $c_1\varepsilon$, where $c_1$ depends on $L$, $\eta$, $\gamma$, and an upper bound $\beta$ on the norm of the Hessian of a quadratic term used in the model $m_k$, the latter of which is equal to $L$ if an exact second-order Taylor series model is used.  The reduction in $f$ on successful iterations can be seen as follows.  First, Cauchy decrease (see Lemma~4.3 in Nocedal \& Wright (2006)) and the bound $\alpha_k \leq \tau \|\nabla f(x_k)\|$ imply that
  \bequationNN
    \baligned
      m_k(x_k) - m_k(s_k(\alpha_k))
        &\geq \thalf \|\nabla f(x_k)\| \min\{ \alpha_k, \tfrac{1}{\beta}\|\nabla f(x_k)\| \} \\
        &\geq \thalf (\tfrac{1}{\tau}) \alpha_k \min\{ \alpha_k, \tfrac{1}{\beta\tau} \alpha_k\} \\
        &=    \thalf (\tfrac{1}{\tau} \min\{ 1, \tfrac{1}{\beta\tau}\}) \alpha_k^2,
    \ealigned
  \eequationNN
  from which it follows that for a successful iteration
  \bequationNN
    f(x_k) - f(x_{k+1}) \geq \thalf \eta (\tfrac{1}{\tau} \min\{ 1, \tfrac{1}{\beta\tau}\}) \alpha_k^2.
  \eequationNN
  \item Analysis in \S\ref{sec.det_ls}.  The lower bound for the reduction in $f$ can be seen in various settings.  For example, if $d_k = -M_k\nabla f(x_k)$ for all $k \in \N{}$, where $\{M_k\}$ is a sequence of real symmetric matrices with eigenvalues uniformly bounded in a positive interval $[\kappa_1,\kappa_2]$, then one finds that
  \bequationNN
    -\eta \nabla f(x_k)^Ts_k(\alpha_k) = -\eta \alpha_k \nabla f(x_k)^Td_k = \eta \alpha_k \nabla f(x_k)^TM_k\nabla f(x_k) \geq \eta \alpha_k \kappa_1 \|\nabla f(x_k)\|^2.
  \eequationNN
  Similarly, the lower bound for $\{\alpha_k\}$ can be seen in various settings.  For example, in the same setting as above, one has from Lipschitz continuity of the gradient that
  \bequationNN
    \baligned
      f(x_k + s_k(\alpha_k)) - f(x_k)
        &\leq \nabla f(x_k)^Ts_k(\alpha_k) + \thalf L \|s_k(\alpha_k)\|^2 \\
        &=    -\alpha_k \nabla f(x_k)^TM_k\nabla f(x_k) + \thalf L \alpha_k^2 \|M_k\nabla f(x_k)\|^2 \\
        &\leq -\alpha_k \kappa_1 \|\nabla f(x_k)\|^2 + \thalf L \kappa_2^2 \alpha_k^2 \|\nabla f(x_k)\|^2.
    \ealigned
  \eequationNN
  On the other hand, if $\alpha_k$ fails to satisfy the Armijo condition, then
  \bequationNN
    \baligned
      f(x_k + s_k(\alpha_k)) - f(x_k) 
        &> \eta \nabla f(x_k)^Ts_k(\alpha_k) \\
        &= - \eta \alpha_k \nabla f(x_k)^TM_k\nabla f(x_k) \\
        &\geq -\eta \alpha_k \kappa_1 \|\nabla f(x_k)\|^2.
    \ealigned
  \eequationNN
  Combined, this shows that
  \bequationNN
    \alpha_k \geq \tfrac{2(1 - \eta)\kappa_1}{L\kappa_2^2},
  \eequationNN
  from which it follows by the structure of the algorithm that $\{\alpha_k\} \geq 2\gamma^{-1}(1-\eta)\kappa_1/(L\kappa_2^2)$.  Finally, the upper bound on the norm of the gradient after a successful step follows since
  \bequationNN
    \baligned
      \|\nabla f(x_{k+1})\|
        &= \|\nabla f(x_{k+1}) - \nabla f(x_k) + \nabla f(x_k)\| \\
        &\leq L\|x_{k+1} - x_k\| + \|\nabla f(x_k)\| \\
        &\leq L\alpha_k\beta\|\nabla f(x_k)\| + \|\nabla f(x_k)\|.
    \ealigned
  \eequationNN
  \item Analysis in \S\ref{sec:extras}.  The desired results in convex and strongly convex settings can be derived using similar techniques as in \cite[\S3.3--\S3.4]{cartis2018global}.  First, suppose~$f$ is convex, has at least one global minimizer (call it $x_*$), and has bounded level sets in the sense that, for some $D \in (0,\infty)$,
  \bequationNN
    \|x - x_*\| \leq D\ \ \text{for all}\ \ x \in \R{n}\ \ \text{with}\ \ f(x) \leq f(x_0).
  \eequationNN
  Convexity of $f$ implies for all $k \in \N{}$ that
  \bequationNN
    f_* - f(x_k) \geq \nabla f(x_k)^T(x_* - x_k) \geq -D\|\nabla f(x_k)\|.
  \eequationNN
  If iteration $k$ is successful, then the above implies that
  \bequationNN
    \baligned
      (f(x_k) - f_*) - (f(x_{k+1}) - f_*)
        &= f(x_k) - f(x_{k+1}) \\
        &\geq \eta c_3 \alpha_k \|\nabla f(x_k)\|^2 \\
        &\geq \tfrac{\eta c_3}{D^2} \alpha_k (f(x_k) - f_*)^2.
    \ealigned
  \eequationNN
  Dividing by $(f(x_k) - f_*)(f(x_{k+1}) - f_*) > 0$ and noting that $f(x_k) > f(x_{k+1})$, one finds that
  \bequationNN
    \tfrac{1}{f(x_{k+1}) - f_*} - \tfrac{1}{f(x_k) - f_*} \geq \tfrac{\eta c_3}{D^2} \alpha_k \tfrac{f(x_k) - f_*}{f(x_{k+1}) - f_*} \geq \tfrac{\eta c_3}{D^2} \underline\alpha,
  \eequationNN
  showing, as desired, that $(f(x_{k+1}) - f_*)^{-1} - (f(x_k) - f_*)^{-1}$ is uniformly bounded below by a positive constant over all successful iterations.
  
  If $f$ is strongly convex, then, as is well known, one has for some $c \in (0,\infty)$ that
  \bequationNN
    \|\nabla f(x)\|^2 \geq c(f(x) - f_*)\ \ \text{for all}\ \ x \in \R{n}.
  \eequationNN
  This shows that if iteration $k$ is successful, then (similar to above)
  \bequationNN
    (f(x_k) - f_*) - (f(x_{k+1}) - f_*) \geq \eta c_3 c \alpha_k (f(x_k) - f_*),
  \eequationNN
  which implies that
  \bequationNN
    f(x_{k+1}) - f_* \leq (1 - \eta c_3 c \alpha_k) (f(x_k) - f_*).
  \eequationNN
  By the definition of $f_*$, it is clear from this inequality that for this successful iteration one must have $\alpha_k \leq 1/(\eta c_3 c)$, meaning $(1 - \eta c_3 c \alpha_k) \in [0,1)$.
  Taking logs, this implies that
  \bequationNN
    \log(f(x_{k+1}) - f_*) \leq \log(1 - \eta c_3 \alpha_k) + \log(f(x_k) - f_*),
  \eequationNN
  which after rearrangement yields
  \bequationNN
    -\log(f(x_{k+1}) - f_*) \geq -\log(1 - \eta c_3 \alpha_k) - \log(f(x_k) - f_*).
  \eequationNN
  As  desired, this shows that $\log((f(x_{k+1}) - f_*)^{-1}) - \log((f(x_k) - f_*)^{-1})$ is uniformly bounded below by a positive constant over all successful iterations.
\eitemize

\fi

\end{document}

%% file: det_good.tex
\begin{tikzpicture}[scale=2,declare function={f(\x)=(1/5)*\x*\x-2*\x+5;g(\x)=(2/5)*\x-2;m(\x)=f(2.0)+g(2.0)*(\x-2.0);}]
  \coordinate (ll) at (0.0, 0.0);
  \coordinate (lr) at (4.0, 0.0);
  \coordinate (ur) at (4.0, 4.0);
  \coordinate (ul) at (0.0, 4.0);
  \draw[black,thick] (ll) -- (lr) -- (ur) -- (ul) -- (ll);
  \draw[red,very thick,smooth,domain=0.5279:4.0] plot({\x},{f(\x)});
  \draw[blue,very thick,smooth,domain=0.17:3.5] plot({\x},{m(\x)});
  \coordinate (xk) at (2.0,0.0);
  \coordinate (xksk) at (3.0,0.0);
  \coordinate [label={\huge $x_k$}] (xklabel) at (2.0,-0.5);
  \coordinate [label={\huge $x_k\!+\!s_k$}] (xksklabel) at (3.0,-0.5);
  \draw[red,very thick] (2.0,{f(2.0)}) circle (0.5mm);
  \draw[red,very thick] (3.0,{f(3.0)}) circle (0.5mm);
  \draw[blue,very thick] (2.0,{m(2.0)}) circle (0.5mm);
  \draw[blue,very thick] (3.0,{m(3.0)}) circle (0.5mm);
  \coordinate (xkskm) at (1.0,0.0);
  \coordinate (xkskm2) at (1.0,4.0);
  \coordinate (xk2) at (2.0,4.0);
  \coordinate (xksk2) at (3.0,4.0);
  \draw[gray,thick] (xkskm) -- (xkskm2);
  \draw[gray,thick,dashed] (xk) -- (xk2);
  \draw[gray,thick] (xksk) -- (xksk2);
\end{tikzpicture}

%% file: det_bad.tex
\begin{tikzpicture}[scale=2,declare function={f(\x)=(1/4)*\x*\x-1.2*\x+3.44;g(\x)=(1/2)*\x-1.2;m(\x)=f(2.0)+0.72*g(2.0)*(\x-2.0);}]
  \coordinate (ll) at (0.0, 0.0);
  \coordinate (lr) at (4.0, 0.0);
  \coordinate (ur) at (4.0, 4.0);
  \coordinate (ul) at (0.0, 4.0);
  \draw[black,thick] (ll) -- (lr) -- (ur) -- (ul) -- (ll);
  \draw[red,very thick,smooth,domain=0.0:4.0] plot({\x},{f(\x)});
  \draw[blue,very thick,smooth,domain=0.0:4.0] plot({\x},{m(\x)});
  \coordinate (xk) at (2.0,0.0);
  \coordinate (xksk) at (3.0,0.0);
  \coordinate [label={\huge $x_k$}] (xklabel) at (2.0,-0.5);
  \coordinate [label={\huge $x_k\!+\!s_k$}] (xksklabel) at (3.0,-0.5);
  \draw[red,very thick] (2.0,{f(2.0)}) circle (0.5mm);
  \draw[red,very thick] (3.0,{f(3.0)}) circle (0.5mm);
  \draw[blue,very thick] (2.0,{m(2.0)}) circle (0.5mm);
  \draw[blue,very thick] (3.0,{m(3.0)}) circle (0.5mm);
  \coordinate (xkskm) at (1.0,0.0);
  \coordinate (xkskm2) at (1.0,4.0);
  \coordinate (xk2) at (2.0,4.0);
  \coordinate (xksk2) at (3.0,4.0);
  \draw[gray,thick] (xkskm) -- (xkskm2);
  \draw[gray,thick,dashed] (xk) -- (xk2);
  \draw[gray,thick] (xksk) -- (xksk2);
\end{tikzpicture}

%% file: good_good.tex
\begin{tikzpicture}[scale=2,declare function={f(\x)=(1/5)*\x*\x-2*\x+5;g(\x)=(2/5)*\x-2;m(\x)=f(2.0)+0.72*g(2.0)*(\x-2.0);}]
  \coordinate (ll) at (0.0, 0.0);
  \coordinate (lr) at (4.0, 0.0);
  \coordinate (ur) at (4.0, 4.0);
  \coordinate (ul) at (0.0, 4.0);
  \draw[black,thick] (ll) -- (lr) -- (ur) -- (ul) -- (ll);
  \draw[red,very thick,smooth,domain=0.5279:4.0] plot({\x},{f(\x)});
  \draw[blue,very thick,smooth,domain=0.0:4.0] plot({\x},{m(\x)});
  \coordinate (xk) at (2.0,0.0);
  \coordinate (xksk) at (3.0,0.0);
  \coordinate [label={\huge $x_k$}] (xklabel) at (2.0,-0.5);
  \coordinate [label={\huge $x_k\!+\!s_k$}] (xksklabel) at (3.0,-0.5);
  \draw[red,very thick] (2.0,{f(2.0)+0.2}) circle (0.5mm);
  \draw[red,very thick] (3.0,{f(3.0)-0.1}) circle (0.5mm);
  \draw[blue,very thick] (2.0,{m(2.0)}) circle (0.5mm);
  \draw[blue,very thick] (3.0,{m(3.0)}) circle (0.5mm);
  \coordinate (xkskm) at (1.0,0.0);
  \coordinate (xkskm2) at (1.0,4.0);
  \coordinate (xk2) at (2.0,4.0);
  \coordinate (xksk2) at (3.0,4.0);
  \draw[gray,thick] (xkskm) -- (xkskm2);
  \draw[gray,thick,dashed] (xk) -- (xk2);
  \draw[gray,thick] (xksk) -- (xksk2);
\end{tikzpicture}

%% file: good_bad.tex
\begin{tikzpicture}[scale=2,declare function={f(\x)=(1/5)*\x*\x-2*\x+5;g(\x)=(2/5)*\x-2;m(\x)=f(2.0)+0.72*g(2.0)*(\x-2.0);}]
  \coordinate (ll) at (0.0, 0.0);
  \coordinate (lr) at (4.0, 0.0);
  \coordinate (ur) at (4.0, 4.0);
  \coordinate (ul) at (0.0, 4.0);
  \draw[black,thick] (ll) -- (lr) -- (ur) -- (ul) -- (ll);
  \draw[red,very thick,smooth,domain=0.5279:4.0] plot({\x},{f(\x)});
  \draw[blue,very thick,smooth,domain=0.0:4.0] plot({\x},{m(\x)});
  \coordinate (xk) at (2.0,0.0);
  \coordinate (xksk) at (3.0,0.0);
  \coordinate [label={\huge $x_k$}] (xklabel) at (2.0,-0.5);
  \coordinate [label={\huge $x_k\!+\!s_k$}] (xksklabel) at (3.0,-0.5);
  \draw[red,very thick] (2.0,{f(2.0)+0.2}) circle (0.5mm);
  \draw[red,very thick] (3.0,{f(3.0)+1.5}) circle (0.5mm);
  \draw[blue,very thick] (2.0,{m(2.0)}) circle (0.5mm);
  \draw[blue,very thick] (3.0,{m(3.0)}) circle (0.5mm);
  \coordinate (xkskm) at (1.0,0.0);
  \coordinate (xkskm2) at (1.0,4.0);
  \coordinate (xk2) at (2.0,4.0);
  \coordinate (xksk2) at (3.0,4.0);
  \draw[gray,thick] (xkskm) -- (xkskm2);
  \draw[gray,thick,dashed] (xk) -- (xk2);
  \draw[gray,thick] (xksk) -- (xksk2);
\end{tikzpicture}

%% file: bad_good.tex
\begin{tikzpicture}[scale=2,declare function={f(\x)=(1/5)*\x*\x-2*\x+5;g(\x)=(2/5)*\x-2;m(\x)=f(2.0)-0.65*g(2.0)*(\x-2.0);}]
  \coordinate (ll) at (0.0, 0.0);
  \coordinate (lr) at (4.0, 0.0);
  \coordinate (ur) at (4.0, 4.0);
  \coordinate (ul) at (0.0, 4.0);
  \draw[black,thick] (ll) -- (lr) -- (ur) -- (ul) -- (ll);
  \draw[red,very thick,smooth,domain=0.5279:4.0] plot({\x},{f(\x)});
  \draw[blue,very thick,smooth,domain=0.0:4.0] plot({\x},{m(\x)});
  \coordinate (xk) at (2.0,0.0);
  \coordinate (xksk) at (1.0,0.0);
  \coordinate [label={\huge $x_k$}] (xklabel) at (2.0,-0.5);
  \coordinate [label={\huge $x_k\!+\!s_k$}] (xksklabel) at (1.0,-0.5);
  \draw[red,very thick] (2.0,{f(2.0)+0.2}) circle (0.5mm);
  \draw[red,very thick] (1.0,{f(1.0)-0.15}) circle (0.5mm);
  \draw[blue,very thick] (2.0,{m(2.0)}) circle (0.5mm);
  \draw[blue,very thick] (1.0,{m(1.0)}) circle (0.5mm);
  \coordinate (xkskm) at (3.0,0.0);
  \coordinate (xkskm2) at (3.0,4.0);
  \coordinate (xk2) at (2.0,4.0);
  \coordinate (xksk2) at (1.0,4.0);
  \draw[gray,thick] (xkskm) -- (xkskm2);
  \draw[gray,thick,dashed] (xk) -- (xk2);
  \draw[gray,thick] (xksk) -- (xksk2);
\end{tikzpicture}

%% file: bad_bad.tex
\begin{tikzpicture}[scale=2,declare function={f(\x)=(1/5)*\x*\x-2*\x+5;g(\x)=(2/5)*\x-2;m(\x)=f(2.0)-0.65*g(2.0)*(\x-2.0);}]
  \coordinate (ll) at (0.0, 0.0);
  \coordinate (lr) at (4.0, 0.0);
  \coordinate (ur) at (4.0, 4.0);
  \coordinate (ul) at (0.0, 4.0);
  \draw[black,thick] (ll) -- (lr) -- (ur) -- (ul) -- (ll);
  \draw[red,very thick,smooth,domain=0.5279:4.0] plot({\x},{f(\x)});
  \draw[blue,very thick,smooth,domain=0.0:4.0] plot({\x},{m(\x)});
  \coordinate (xk) at (2.0,0.0);
  \coordinate (xksk) at (1.0,0.0);
  \coordinate [label={\huge $x_k$}] (xklabel) at (2.0,-0.5);
  \coordinate [label={\huge $x_k\!+\!s_k$}] (xksklabel) at (1.0,-0.5);
  \draw[red,very thick] (2.0,{f(2.0)+0.2}) circle (0.5mm);
  \draw[red,very thick] (1.0,{f(1.0)-1.8}) circle (0.5mm);
  \draw[blue,very thick] (2.0,{m(2.0)}) circle (0.5mm);
  \draw[blue,very thick] (1.0,{m(1.0)}) circle (0.5mm);
  \coordinate (xkskm) at (3.0,0.0);
  \coordinate (xkskm2) at (3.0,4.0);
  \coordinate (xk2) at (2.0,4.0);
  \coordinate (xksk2) at (1.0,4.0);
  \draw[gray,thick] (xkskm) -- (xkskm2);
  \draw[gray,thick,dashed] (xk) -- (xk2);
  \draw[gray,thick] (xksk) -- (xksk2);
\end{tikzpicture}